\documentclass[numbers=noenddot]{scrartcl}
\usepackage [latin1]{inputenc}
\usepackage{amsthm,amsmath,amssymb}
\usepackage[colorlinks=true,citecolor=blue,linkcolor=blue,urlcolor=blue,
  pdfauthor={Thomas Apel, Katharina Lorenz, Serge Nicaise},
  pdfstartview={Fit}]{hyperref}

\theoremstyle{plain}

\theoremstyle{definition}

\theoremstyle{remark}

\usepackage[dvipsnames]{xcolor}
\usepackage[ulem=normalem,draft,authormarkup=none,deletedmarkup=xout]{changes}
\definechangesauthor[name=Katharina, color=ForestGreen]{KL}
\definechangesauthor[name=Thomas, color=NavyBlue]{TA}
\definechangesauthor[name=Serge, color=Magenta]{SN}

\newcommand{\R}{\mathbb{R}}

\numberwithin{equation}{section}

\newcommand{\norm}[1]{\lVert#1\rVert}

\begin{document}
\title{Weak and very weak solutions of the Laplace equation and the
  Stokes system with prescribed regularity} \author{Thomas
  Apel\thanks{\texttt{thomas.apel@unibw.de}, Universit\"at der
    Bundeswehr M\"unchen, Institute of Mathematics and Com\-pu\-ter-Based
    Simulation, 85577 Neubiberg, Germany.} \and Katharina
  Lorenz\thanks{\texttt{katharina.lorenz@unibw.de}, Universit\"at der
    Bundeswehr M\"unchen, Institute of Mathematics and Computer-Based
    Simulation, 85577 Neubiberg, Germany.  } \and Serge
  Nicaise\thanks{\texttt{serge.nicaise@uphf.fr}, Universit\'e
    Polytechnique Hauts-de-France, C\'ERAMATHS/DMATHS and FR CNRS 2037, F-59313 -
    Valenciennes Cedex 9, France.}} 
\maketitle

\begin{small}
  \setlength{\parindent}{0pt}\setlength{\parskip}{1ex}
  \textbf{Abstract:} To verify theoretical results it is sometimes
  important to use a numerical example where the solution has a
  particular regularity. The paper describes one approach to construct
  such examples. It is based on the regularity theory for elliptic
  boundary value problems. 


\end{small}

\section{Introduction}

The design of numerical examples is often a delicate task. The test
problem should confirm the theory but the solution should not be more
regular than assumed in the theory. We formulate here two types of
examples for the Laplace and the Stokes problems. In Subsections
\ref{sec:1.1} and \ref{sec:2.1} we formulate families of solutions of
the homogeneous differential equations which solve problems determined
by the non-homogeneous boundary datum. The solutions are smooth except
in the vicinity of one boundary point, and a parameter controls their
regularity. In Subsections \ref{sec:1.2} and \ref{sec:2.2} the
parameters are restricted to achieve homogeneous boundary conditions,
at least in the vicinity of the singular boundary point.  The solution
is then determined by the boundary datum away from the singularity or,
by using a cut-off function, by a smooth right hand side of the
differential equation.

In the case of the Laplace operator, these examples are widely known
but not so much for the Stokes problem. In order to explain the ideas
we start, however, with the Laplace equation. We assume that the
domain $\Omega$ is two-dimensional and polygonal, and comment on the
three-dimensional case in Subsections \ref{rem:3dLaplace} and
\ref{rem:3dstokes}.

\section{The Laplace equation}
In Subsections \ref{sec:1.1} to \ref{sec:test} we assume that  $\Omega$ is a bounded polygon with a Lipschitz boundary.

\subsection{\label{sec:1.1}Fundamental solutions}
Let $(r,\theta)$ be polar coordinates centered in the corner of
$\Omega\subset\R^2$ with maximal interior angle $\omega$ such that the
edges of this corner are described by $\theta=0$ and $\theta=\omega$.
Then the function
\begin{align*}
  u(r,\theta)=r^\lambda\Phi(\theta)  
\end{align*}
solves the Laplace equation 
\begin{align}\label{eq:Laplace}
  -\Delta u&=0\quad\text{in }\Omega
\end{align}
iff $\Phi''+\lambda^2\Phi=0$, i.\,e., iff
\begin{align*}
  \Phi(\theta)=\begin{cases}
    c_1\cos\lambda\theta+c_2\sin\lambda\theta & \text{if }\lambda\not=0, \\
    c_1+c_2\theta & \text{if }\lambda=0,
  \end{cases}  
\end{align*}
such that the solution is
\begin{align}\label{eq:stern}
  u(r,\theta)=\begin{cases}
    r^\lambda(c_1\cos\lambda\theta+c_2\sin\lambda\theta) & 
    \text{if }\lambda\not=0, \\
    c_1+c_2\theta & \text{if }\lambda=0.
  \end{cases}   
\end{align}
This solution satisfies
\begin{align*}
  u\in H^s(\Omega) \quad\forall s<1+\lambda
\end{align*}
such that the choice of $\lambda$ can be used to develop test examples of the form
\begin{align}
  -\Delta u&=0\quad\text{in }\Omega, \label{eq:Laplace1}\\
  u&=g\quad\text{on }\Gamma:=\partial\Omega, \label{eq:DirichletRB}
\end{align}
or with other types of boundary conditions, with solutions of
prescribed regularity.

The case $\lambda=0$ is of interest when the boundary datum
should have a jump at $r=0$.

\subsection{\label{sec:1.2}Boundary conditions}
The freedom of choosing $\lambda$, $c_1$, and $c_2$ can also be used
to satisfy homogeneous boundary conditions, e.\,g.,
\begin{align*}
  c_1&=0,&\lambda&=k\tfrac\pi\omega,& k&=1,2, \ldots, && 
  \quad\text{leads to } u(r,0)=u(r,\omega)=0, \\
  c_2&=0,&\lambda&=k\tfrac\pi\omega,& k&=0,1,2, \ldots, && 
  \quad\text{leads to } \partial_n u(r,0)=\partial_n u(r,\omega)=0, \\
  c_1&=0,&\lambda&=(k-\tfrac12)\tfrac\pi\omega,& k&=1,2, \ldots, && 
  \quad\text{leads to } u(r,0)=\partial_n u(r,\omega)=0.
\end{align*}
The remaining constant is still free. The function does not vanish at
boundary parts which are not adjacent to the point with $r=0$. If
homogeneous boundary conditions should be satisfied on the whole
boundary of the domain, the functions could be multiplied by a smooth
cut-off function $\eta:\R_+\to\R$ with
\begin{align} \label{eq:cut-off}
  \eta(r)=\begin{cases}
    1 & 
    \text{if }r<r_0, \\
    0 & \text{if }r>r_1>r_0,
  \end{cases}   
\end{align}
and appropriately chosen positive constants $r_0,r_1\in\R$. Then $f=-\Delta u$
is zero only for $r\le r_0$ and $r\ge r_1$ but smooth in the region
$r_0<r<r_1$.

Note that for a given domain $\Omega$ the regularity can be influenced
only in discrete steps with this approach. But by adjusting the
interior angle of $\Omega$ at $r=0$, any desired regularity is adjustable.

\subsection{\label{sec:1.3}Weak and very weak solutions}

For $\lambda>0$, the function $u$ from \eqref{eq:stern} is a weak
solution of \eqref{eq:Laplace}, in the sense that it belongs to
$H^1(\Omega)$ and satisfies
\begin{align*}
  (\nabla u,\nabla v)&=0 \quad\forall v\in H^1_0(\Omega),
\end{align*}
and for $-\min(1,\xi)<\lambda\le0$, $\xi:=\frac\pi\omega$, it is a
very weak solution, in the sense that it belongs to $L^2(\Omega)$
and satisfies
\begin{align*}
  (u,\Delta v)=\langle u,\partial_nv\rangle_\Gamma \quad\forall
  v\in V=\{v\in H^1_0(\Omega)\colon\Delta v\in L^2(\Omega)\},
\end{align*}
see also \cite{ApelNicaisePfefferer:16}. 

If $\lambda\le-1$, then the function $u$ is not a very weak solution
since $u\not\in L^2(\Omega)$ such that $(u,\Delta v)$ is not well
defined for all $v\in V$.

If the domain $\Omega$ is non-convex and $-1<\lambda<-\xi$, then for
$v=\eta(r) r^\xi \sin\xi\theta \in V$ and for
$\Gamma_0=\{x\in\Gamma:r(x)<r_0\}$ the dual bracket $\langle
u,\partial_nv\rangle_{\Gamma_0}$ is meaningless because the product
$u\partial_n v$ is not Lebesgue integrable. Indeed, on $\theta=0$, we
have
\[
  |u\partial_n v|=|\Phi(0)|\,\eta(r)\, r^{\lambda+\xi-1},
\]
which is not integrable except if $\Phi(0)=0$. Hence $u\partial_n v$
is integrable on $\Gamma_0$ if and only if $\Phi(0)=\Phi(\omega)=0$,
which is not possible since $\lambda\ne -k\xi$ for all $k=1,2,
\ldots$.

Let us finally consider the limit case
$u^*(r,\theta)=r^{-\xi}\sin\xi\theta$. Note first that $\Delta u^*=0$
in $\Omega$ and $u^*\in H^t(\Omega)$ for all $t<1-\xi$.
%
%
However, we can show  
that
\[
(u^*,\Delta v)-\langle u^*,\partial_nv\rangle_\Gamma\ne 0
\]
for $v=\eta(r) r^\xi \sin\xi\theta \in V$. Indeed,
for any $\varepsilon>0$, by setting
$
\Omega_\varepsilon=\Omega\setminus \bar B(0,\varepsilon)$, we may write
\[
(u^*,\Delta v)=\lim_{\varepsilon\to 0+}
\int_{\Omega_\varepsilon} u\Delta v.
\]
Now since $u^*$ and $v$ are smooth in $\Omega_\varepsilon$, we can use the Green formula to show that
\[
(u^*,\Delta v)  
=\lim_{\varepsilon\to 0+}\int_{\partial \Omega_\varepsilon} 
(u^*\partial_n v-v \partial_n u^*).
\]
Splitting the integral in $\partial \Omega_\varepsilon$ into the 
integral into the boundary of the disc $B(0,\varepsilon)$
and the remainder part, and using the form of $u$ and $v$, we get
\[
(u^*,\Delta v)-\langle u^*,\partial_nv\rangle_\Gamma
=\lim_{\varepsilon\to 0+}\int_{0}^\omega (u^*\partial_n v-v \partial_n u^*)\varepsilon\,\mathrm{d}\theta
=-2\xi \int_{0}^\omega \sin^2(\xi \theta)\,\mathrm{d}\theta=
-\pi,
\]
which proves the assertion. Alternatively, one could show that the
problem \eqref{eq:DirichletRB} with $g=u^*$ on $\Gamma$ has a weak
solution $u\not=u^*$ which also means that $u^*$ is not a weak or very
weak solution, see Subsection \ref{sec:test}. 

\subsection{\label{sec:test}Numerical test for the limit case}

\begin{figure}[ht] %
  \centering
  \includegraphics[width=0.6\textwidth]{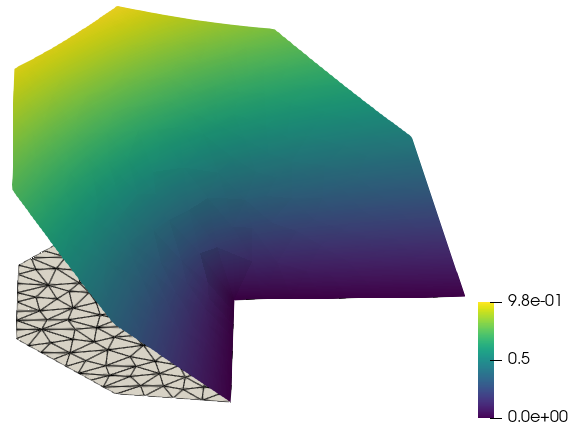}
  \caption{\label{fig:test}Illustration of $u^*$}
\end{figure}

The numerical test is carried out using the method described in \cite[Section III.B]{ApelNicaisePfefferer:16} which proved to be able to approximate very weak solutions, even with pole, see \cite[Section IV]{ApelNicaisePfefferer:16}.
For the illustration here, we consider a polygonal domain \( \Omega \subset \R^2 \) with an interior angle \( \omega = \frac{3}{2} \pi \). Accordingly, we set \( \xi = \frac{2}{3} \), and the boundary datum $g$ is given by  
\[
g(r,\theta) = r^{-\frac{2}{3}} \sin(\tfrac{2}{3} \theta).
\]
It fulfills $g(r,0) = g(r,\frac{3}{2} \pi) = 0$ and is not zero at the boundary parts which are not adjacent to $(0,0)$. 
However, the obtained solution $u(r,\theta)$ is not $r^{-\frac{2}{3}} \sin(\tfrac{2}{3} \theta)$. In particular, it has no pole,
see the illustration in Figure \ref{fig:test}.
The explanation is that the boundary datum is piecewise smooth and continuous, hence at least in $H^{\frac12}(\Gamma)$ such that a unique weak solution $u\in  H^1(\Omega)$ exists.
In conclusion, this numerical test confirms that we are in the limit case, where $r^{-\frac{2}{3}} \sin(\tfrac{2}{3} \theta)\in L^2(\Omega)$ is not a very weak solution although is harmonic and solves the Laplace equation.

\subsection{Three-dimensional case \label{rem:3dLaplace}}
The boundary of a polyhedral domain has edges and corners. Near an
edge, the typical behavior of a solution of the Poisson equation is
\begin{align*}
cr^\lambda\Phi(\theta)\in H^s(\Omega),\quad s<1+\lambda,
\end{align*}
with some function $c$, where we used polar
coordinates $(r,\theta)$ perpendicular to the edge. Boundary
conditions can be satisfied as in the two-dimensional case described
in Subsection \ref{sec:1.2}. 

A second family of fundamental solutions has to be considered near
corners. Using spherical coordinates $(R,\vartheta,\theta)$, they can
be described by 
\begin{align*}
  R^\nu\Phi_c(\vartheta,\theta)\in H^s(\Omega),\quad
  s<\tfrac32+\nu
\end{align*}
if $\Phi_c$ is a smooth enough solution of
$\Delta'\Phi_c+\nu(\nu+1)\Phi_c=0$ in the intersection of $\Omega$
with a sphere centered at the corner, $\Delta'$ being the
Laplace--Beltrami operator.  Boundary conditions can be satisfied when
$\Phi_c$ is an eigenfunction of the Laplace--Beltrami operator with
corresponding boundary conditions defined on the intersection of
$\Omega$ with a sphere centered at the corner; however this can be
done analytically only in very special cases \cite[\S 18D,
\S18E]{dauge:88}. Alternatively, the pair $(\nu,\Phi_c)$ can be
approximated numerically, see \cite{pester:06a} and the literature
cited therein.

\section{The Stokes system}
In Subsections \ref{sec:2.1} to \ref{sec:2.3} we assume that  $\Omega$ is a bounded polygon with a Lipschitz boundary.

\subsection{\label{sec:2.1}Fundamental solutions}

We describe now solutions of the Stokes system of the form
\begin{align}\label{eq:ansatz}
  u(r,\theta)=r^\lambda U(\theta), \quad p(r,\theta)=r^{\lambda-1} P(\theta),
\end{align}
where again $(r,\theta)$ with $r\in\mathbb{R}_+$, $\theta\in(0,\omega)$ are polar
coordinates centered in the corner with maximal interior angle, such that
\begin{align*}
  -\Delta u+\nabla p &= 0 \quad\text{in }\Omega\subset\R^2, \\
  \nabla\cdot u&=0\quad \text{in }\Omega.
\end{align*}
It is now important to write $u$ in polar components,
$u=u_re_r+u_\theta e_\theta$ with $e_r=\cos\theta\,e_1
+\sin\theta\,e_2$, $e_\theta=-\sin\theta\,e_1+\cos\theta\,e_2$, such
that $\norm{e_r}=\norm{e_\theta}=1$. Note that the 
components $u_r$ and $u_\theta$ are related to the Cartesian
components of $u$ via a rotation computed by
\begin{align*}
  \begin{pmatrix} u_1\\u_2\end{pmatrix} = 
  \begin{pmatrix}\cos\theta&-\sin\theta\\\sin\theta&\cos\theta\end{pmatrix}
  \begin{pmatrix}u_r\\u_\theta\end{pmatrix}.
\end{align*}

We follow the derivation in \cite[Sect. 5.1]{KozlovMazyaRossmann:01}
and \cite[Sect. 9.3.2]{MazyaRossmann:10}.  Inserting the ansatz
\eqref{eq:ansatz} into the equations, the variable $r$ cancels out,
and it remains a system of ordinary differential equations for
$U$ and $P$. The four solutions of this system  are for
$\theta\not=0$  
\begin{align*}
  U_r^{(1)} &= \cos(1+\lambda)\theta,  & U_\theta^{(1)}&=-\sin(1+\lambda)\theta, &  P^{(1)} &=0, \\
  U_r^{(2)} &= \sin(1+\lambda)\theta,  & U_\theta^{(2)}&=\cos(1+\lambda)\theta,  &  P^{(2)} &=0, \\
  U_r^{(3)} &= (1-\lambda)\cos(1-\lambda)\theta, & U_\theta^{(3)}&=-(1+\lambda)\sin(1-\lambda)\theta, &  P^{(3)} &=-4\lambda\cos(1-\lambda)\theta, \\
  U_r^{(4)} &= (1-\lambda)\sin(1-\lambda)\theta, & U_\theta^{(4)}&= (1+\lambda)\cos(1-\lambda)\theta, &  P^{(4)} &=-4\lambda\sin(1-\lambda)\theta.
\end{align*}
For the convenience of the reader, we convert the velocities to
Cartesian components in the Appendix.
In the case $\lambda=0$, the functions $(U^{(i)},P^{(i)})$, $i=1,2$, remain the same (with simplifications), but the others have to be replaced by
\begin{align*}
  U_r^{(3)} &= -\cos\theta+2\theta\sin\theta & U_\theta^{(3)}&=-\sin\theta+2\theta\cos\theta, & P^{(3)} &=-4\cos\theta, \\
  U_r^{(4)} &= -\sin\theta-2\theta\cos\theta & U_\theta^{(4)}&=\cos\theta+2\theta\sin\theta,  & P^{(4)} &=-4\sin\theta.
\end{align*}
The solutions $u$, $p$, from \eqref{eq:ansatz} can now be concluded with
\begin{align*}
  U(\theta)=\sum_{i=1}^4 c_iU^{(i)}(\theta), \quad
  P(\theta)= \sum_{i=1}^4 c_iP^{(i)}(\theta),
\end{align*}
with arbitrary $\lambda$ and arbitrary coefficients $c_i$. They satisfy 
\begin{align*}
  u\in H^s(\Omega),\quad p\in H^{s-1}(\Omega) \quad\forall s<1+\lambda,
\end{align*}
and the parameter $\lambda$ can be chosen such that the test example
has the desired regularity. 

\subsection{\label{sec:2.2}Boundary conditions}

As in Subsection \ref{sec:1.2}, the coefficients $c_i$ and the
parameter $\lambda$ can be used to satisfy homogeneous boundary
conditions. Two boundary conditions for both $\theta=0$ and
$\theta=\omega$ give a homogeneous linear system of 4 equations which
has a non-trivial solution iff the determinant vanishes. This condition
is used to find again a countable number of values of $\lambda$. Let
us sketch this approach for the case of Dirichlet boundary conditions
and $\lambda\not=0$.

The condition $U(0)=0$ leads to
\begin{align*}
  \begin{pmatrix}U_r(0)\\U_\theta(0)\end{pmatrix}&=
  \begin{pmatrix}c_1+(1-\lambda)c_3\\c_2 + (1+\lambda)c_4\end{pmatrix}=
  \begin{pmatrix}0\\0\end{pmatrix},\qquad\text{i.\,e. } &
  \begin{pmatrix}c_1\\c_2\end{pmatrix}&=
  \begin{pmatrix}-(1-\lambda)c_3\\-(1+\lambda)c_4\end{pmatrix},
\end{align*}
hence
\begin{align*}
  U(\omega)=-(1-\lambda)c_3U^{(1)}(\omega)
  -(1+\lambda)c_4U^{(2)}(\omega)+
  c_3U^{(3)}(\omega)+c_4U^{(4)}(\omega).
\end{align*}
The $2\times 2$ linear system $U(\omega)=0$ for the coefficients $c_3$ and $c_4$ has the determinant 
\begin{align}\label{eq:lambda_k}
  4(\sin^2\lambda\omega-\lambda^2\sin^2\omega)=
  4(\sin\lambda\omega-\lambda\sin\omega)(\sin\lambda\omega+\lambda\sin\omega).
\end{align}
This means that for given angle $\omega$ one gets the corresponding
exponents $\lambda\in\mathbb{C}$ by solving (separately) the two transcendental,
scalar equations $\sin\lambda\omega=\pm\lambda\sin\omega$. All values
$\mathrm{Re}\lambda\in[\frac12,4]$ are given for $\omega_k=k\pi/10$,
$k=4,5,\ldots,20$, in \cite{dauge:89}. 

\subsection{\label{sec:2.3}Weak and very weak solutions}
As in Subsection \ref{sec:1.3}, the pair $(u,p)$ is a weak solution
for $\lambda>0$ and a very weak solution for
$-\min(1,\xi)<\lambda\le0$, where  
\begin{align*}
  \xi=\min\{\mathrm{Re}\lambda>0\colon 
  \lambda\text{ satisfies \eqref{eq:lambda_k}}\}
\end{align*}
in the case of Dirichlet boundary conditions. The weak solution
$(u,p)\in H^1(\Omega)\times L^2_0(\Omega)$ is defined by
\begin{align*}
  (\nabla u,\nabla v)-(\nabla\cdot v,p) &= 0 \quad\forall v\in H^1_0(\Omega), \\
  (\nabla\cdot u,q) &=0 \quad\forall q\in L^2_0(\Omega). 
\end{align*}
The very weak solution $(u,p)\in L^2(\Omega)\times P'$ with $P=\{v\in
H^1(\Omega)\cap L^2_0(\Omega): r^{-1} v\in L^2(\Omega)\}$ is defined
by
\begin{align*}
  (u,-\Delta v+\nabla q)-(\nabla\cdot v,p)&= \langle u,qn-\partial_n
  v\rangle_\Gamma \quad \forall (v,q)\in \mathcal{V}
\end{align*}
where $\mathcal{V}:=\{(v,q)\in H^1_0(\Omega)\times L^2_0(\Omega)\colon
-\Delta v+\nabla q \in L^2(\Omega), \nabla\cdot v\in P\}$, see \cite{ALP:24}.

\subsection{Three-dimensional case \label{rem:3dstokes}}

The three-dimensional case is very similar to the Laplace case
described in Subsection \ref{rem:3dLaplace}. The most interesting
difference is that the function $u$ is a vector function with
different regularities of the components near edges. The components
$u_r$ and $u_\theta$ have the two-dimensional behaviour as described
in Subsections \ref{sec:2.1} and \ref{sec:2.2}.  The component in edge
direction, however, has a behaviour as for the Laplace operator. For
more details, see e.\,g. \cite[Chap. 9]{MazyaRossmann:10}.


\appendix

\section{Velocities in Cartesian coordinates}
  Since most finite element packages work with Cartesian components of
  the vector functions, we convert here the velocities in the
  fundamental solutions of the Stokes system.
\begin{align*}
  \begin{pmatrix}U_1^{(1)}\\U_2^{(1)}\end{pmatrix} &= 
  \begin{pmatrix}\cos\theta&-\sin\theta\\\sin\theta&\cos\theta\end{pmatrix} 
  \begin{pmatrix}\cos(1+\lambda)\theta\\-\sin(1+\lambda)\theta\end{pmatrix} \\
  &=\begin{pmatrix}\cos\theta\cos(1+\lambda)\theta+\sin\theta\sin(1+\lambda)\theta \\
  \sin\theta\cos(1+\lambda)\theta-\cos\theta\sin(1+\lambda)\theta\end{pmatrix} 
  =\begin{pmatrix}\cos\lambda\theta \\ - \sin\lambda\theta\end{pmatrix} \\
  \begin{pmatrix}U_1^{(2)}\\U_2^{(2)}\end{pmatrix} &= 
  \begin{pmatrix}\cos\theta&-\sin\theta\\\sin\theta&\cos\theta\end{pmatrix} 
  \begin{pmatrix}\sin(1+\lambda)\theta\\\cos(1+\lambda)\theta\end{pmatrix} \\
  &=\begin{pmatrix}\cos\theta\sin(1+\lambda)\theta-\sin\theta\cos(1+\lambda)\theta \\
  \sin\theta\sin(1+\lambda)\theta+\cos\theta\cos(1+\lambda)\theta\end{pmatrix} 
  =\begin{pmatrix}\sin\lambda\theta \\ \cos\lambda\theta\end{pmatrix} \\
  \begin{pmatrix}U_1^{(3)}\\U_2^{(3)}\end{pmatrix} &= 
  \begin{pmatrix}\cos\theta&-\sin\theta\\\sin\theta&\cos\theta\end{pmatrix} 
  \begin{pmatrix}(1-\lambda)\cos(1-\lambda)\theta\\-(1+\lambda)\sin(1-\lambda)\theta\end{pmatrix} \\
  &=\begin{pmatrix}(1-\lambda)\cos\theta\cos(1-\lambda)\theta+(1+\lambda)\sin\theta\sin(1-\lambda)\theta \\
  (1-\lambda)\sin\theta\cos(1-\lambda)\theta-(1+\lambda)\cos\theta\sin(1-\lambda)\theta\end{pmatrix} \\ 
  &=\begin{pmatrix}\cos\lambda\theta-\lambda\cos(2-\lambda)\theta \\ \sin\lambda\theta-\lambda\sin(2-\lambda)\theta\end{pmatrix} \\
  \begin{pmatrix}U_1^{(4)}\\U_2^{(4)}\end{pmatrix} &= 
  \begin{pmatrix}\cos\theta&-\sin\theta\\\sin\theta&\cos\theta\end{pmatrix} 
  \begin{pmatrix}(1-\lambda)\sin(1-\lambda)\theta\\(1+\lambda)\cos(1-\lambda)\theta\end{pmatrix} \\
  &=\begin{pmatrix}(1-\lambda)\cos\theta\sin(1-\lambda)\theta-(1+\lambda)\sin\theta\cos(1-\lambda)\theta \\
  (1-\lambda)\sin\theta\sin(1-\lambda)\theta+(1+\lambda)\cos\theta\cos(1+\lambda)\theta\end{pmatrix} \\
  &=\begin{pmatrix}-\sin\lambda\theta-\lambda\sin(2-\lambda)\theta \\ \cos\lambda\theta+\lambda\cos(2-\lambda)\theta\end{pmatrix} 
\end{align*}
In the case $\lambda=0$ we modify $U^{(3)}$ and $U^{(4)}$ to
\begin{align*}
  \begin{pmatrix}U_1^{(3)}\\U_2^{(3)}\end{pmatrix} &= 
  \begin{pmatrix}\cos\theta&-\sin\theta\\\sin\theta&\cos\theta\end{pmatrix} 
  \begin{pmatrix}-\cos\theta+2\theta\sin\theta\\-\sin\theta+2\theta\cos\theta\end{pmatrix} \\
  &=\begin{pmatrix}\cos\theta(-\cos\theta+2\theta\sin\theta)-\sin\theta(-\sin\theta+2\theta\cos\theta) \\
  \sin\theta(-\cos\theta+2\theta\sin\theta)+\cos\theta(-\sin\theta+2\theta\cos\theta)\end{pmatrix} 
  =\begin{pmatrix}-\cos2\theta \\ 2\theta-\sin2\theta\end{pmatrix} \\
  \begin{pmatrix}U_1^{(4)}\\U_2^{(4)}\end{pmatrix} &= 
  \begin{pmatrix}\cos\theta&-\sin\theta\\\sin\theta&\cos\theta\end{pmatrix} 
  \begin{pmatrix}-\sin\theta-2\theta\cos\theta\\\cos\theta+2\theta\sin\theta\end{pmatrix} \\
  &=\begin{pmatrix}\cos\theta(-\sin\theta-2\theta\cos\theta)-\sin\theta(\cos\theta+2\theta\sin\theta) \\
  \sin\theta(-\sin\theta-2\theta\cos\theta)+\cos\theta(\cos\theta+2\theta\sin\theta)\end{pmatrix} 
  =\begin{pmatrix}-2\theta-\sin2\theta \\\cos2\theta\end{pmatrix} 
\end{align*}

\bibliographystyle{alphaurl}
\bibliography{stokes}

\end{document}